\documentclass[reqno,11pt]{amsart}
\usepackage{amsmath, latexsym, amsfonts, amssymb, amsthm, amscd}

\numberwithin{equation}{section}

\newtheorem{theorem}{Theorem}[section]
\newtheorem{lemma}[theorem]{Lemma}

\newtheorem{cor}[theorem]{Corollary}

\newcommand{\ga}{\alpha}

\newcommand{\gga}{\gamma}            

\newcommand{\gep}{\varepsilon}       

\newcommand{\gl}{\lambda}

\newcommand{\go}{\omega}

\newcommand{\cZ}{{\ensuremath{\mathcal Z}} }

\newcommand{\bbE}{{\ensuremath{\mathbb E}} }

\newcommand{\bbP}{{\ensuremath{\mathbb P}} }

\newcommand{\bbR}{{\ensuremath{\mathbb R}} }

\title[The Brownian meander and excursion conditioned to have a fixed
average]
{On the Brownian meander and excursion conditioned to have a fixed
time average}

\author{Lorenzo Zambotti}
\address{Laboratoire de Probabilit{\'e}s et Mod\`eles Al\'eatoires --
Universit{\'e} Pierre et Marie Curie -  Paris 6 -- Case 188, 4 place
Jussieu, 75252 Paris cedex 05, France }
\email{zambotti\@@ccr.jussieu.fr}
\urladdr{http://www.proba.jussieu.fr/\raisebox{0.15ex}{\tiny$\sim$}zambotti/}

\date{}

\begin{document}

\begin{abstract}
We study the density of the time average of the Brownian
meander/excursion over the time interval $[0,1]$. Moreover we give
an expression for the law of the Brownian meander/excursion
conditioned to have a fixed time average.
\end{abstract}

\maketitle

\section{Introduction}

Let $(m_t,t\in[0,1])$ be a standard Brownian meander and
set $\langle m,1\rangle:=\int_0^1 m_r \, dr$, average of $m$.
In this note we answer two questions:
\begin{enumerate}
\item what is the density of the random variable $\langle m,1\rangle$?
\item what is the law of $(m_t,t\in[0,1])$ conditionally on
the value of $\langle m,1\rangle$?
\end{enumerate}
We provide an expression for both objects,
based on suitable Girsanov transformations. We recall that
$m$ is equal in law to a Brownian motion $B$ conditioned to be
non-negative on $[0,1]$. Using this result, we compute first
the law of $B$ conditioned to have a fixed average, then we
write a Girsanov transformation and finally we condition
this process to be non-negative.

We answer similar questions for the normalized Brownian excursion.
We refer to \cite{reyo} for details about the Brownian
meander and the Brownian excursion.

\subsection{The result}

Let $m$ be a Brownian meander on $[0,1]$ and
$B$ a standard Brownian motion such that $\{m,B\}$ are
independent and let $c\geq 0$ be a constant.
We introduce the continuous processes:
\begin{align*}
& u_t \, := \, \left\{ \begin{array}{ll}
{\displaystyle
\frac 1{\sqrt 2}\, m_{2t}, \qquad \quad t\in[0,1/2]
}
\\ \\
{\displaystyle
\frac 1{\sqrt 2}\, m_1 + B_{t-\frac 12}, \quad t\in[1/2,1],}
\end{array} \right.
\\
& U^c_t \, := \,
\left\{ \begin{array}{ll}
{\displaystyle
u_t, \qquad \quad t\in[0,1/2]
}
\\ \\
{\displaystyle u_t
+ \left(12\, t\, (2-t)-9\right)\left(c-\int_0^1 u\right), \quad t\in[1/2,1].}
\end{array} \right.
\end{align*}
Notice that $\int_0^1 U^c_t \, dt = c$.
\begin{theorem}\label{1.1}
For all bounded Borel $\Phi:C([0,1])\mapsto\bbR$ and $f:\bbR\mapsto\bbR$:
\begin{equation}\label{e1.1}
 \bbE\left[\Phi(m) \, f(\langle m,1\rangle)\right] =
 \int_0^\infty \sqrt\frac{24}\pi \,
\bbE\left[\Phi(U^c)\,
e^{-12\left(\int_0^{1/2}(U_r^c+U_{1/2}^c)\,
dr-c\right)^2}\, 1_{\{U^c_t\geq 0, \ \forall t\in[0,1]\}}\right] f(c) \, dc.
\end{equation}
\end{theorem}
\begin{cor}\label{cor1}
If we set for all $c\geq 0$
\[
p_{\langle m,1\rangle}(c) \, := \, \sqrt\frac{24}\pi \,
\bbE\left[e^{-12\left(\int_0^{1/2}(U_r^c+U_{1/2}^c)\, dr-c\right)^2}
\, 1_{\{U^c_t\geq 0, \ \forall t\in[0,1]\}}\right],
\]
then $p_{\langle m,1\rangle}$ is the density of $\langle
m,1\rangle$, i.e.
\[
\bbP(\langle m,1\rangle\in dc)\, = \, p_{\langle m,1\rangle}(c)
\,1_{\{c\geq 0\}} \, dc.
\]
Moreover $p_{\langle m,1\rangle}$ is continuous on $[0,\infty)$.
\end{cor}
\noindent Notice that a.s. $U^c_t\geq 0$ for all $t\in[0,1/2]$,
since a.s. $m\geq 0$: therefore a.s.
\[
\left\{U^c_t\geq 0, \ \forall t\in[0,1]\right\} \, = \,
\left\{U^c_t\geq 0, \ \forall t\in[1/2,1]\right\}.
\]
The probability of this event is positive for all $c>0$ while it is
0 for $c=0$, since $\int_0^1 U^0_t \, dt = 0$. In particular
$p_{\langle m,1\rangle}(0)=0$. Moreover for all $c>0$ we can define
the probability measure on $C([0,1])$ given by
\[
\bbE\left[\Phi(m) \, | \, \langle m,1\rangle =c \right] \, := \,
\frac 1{\cZ_c} \ \bbE\left[\Phi(U^c)\,
e^{-12\left(\int_0^{1/2}(U_r^c+U_{1/2}^c)\, dr-c\right)^2} \,
1_{\{U^c_t\geq 0, \ \forall t\in[0,1]\}}\right],
\]
where $\Phi:C([0,1])\mapsto\bbR$ is bounded Borel and $\cZ_c>0$ is a
normalization factor.
\begin{cor}\label{cor2}
$(\bbP\left[m\in\cdot \, | \, \langle m,1\rangle =c \right], c>0)$
is a regular conditional distribution of $m$ given $\langle
m,1\rangle$, i.e.
\[
\bbP(m\in\cdot\, , \langle m,1\rangle\in dc ) \, = \,
\bbP\left[m\in\cdot \, | \, \langle m,1\rangle =c \right] \
p_{\langle m,1\rangle}(c) \,1_{\{c> 0\}} \, dc.
\]
Moreover $(0,\infty)\ni c\mapsto \bbP\left[m\in\cdot \, | \, \langle
m,1\rangle =c \right]$ is continuous in the weak topology.
\end{cor}
\noindent

In section \ref{exx} below we give similar results for the Brownian
excursion, see Theorem \ref{4.1}.

\medskip
The results of this note are of interest in connection with a class
of Stochastic Partial Differential Equations, studied in \cite{deza}
and \cite{z}, which have the property of conservation of the space
average. For instance, the stochastic equation considered in
\cite{z} admits as invariant measure the Brownian excursion
conditioned to have a fixed average. Moreover, the Brownian meander
conditioned to have a fixed average and the density $p_{\langle
m,1\rangle}$ appear in an infinite-dimensional integration by parts
formula in \cite[Corollary 6.2]{deza}.

\section{An absolute continuity formula}
\label{secabco}

Let $(X_t)_{t\in[0,1]}$ be a continuous centered Gaussian process
with covariance function $q_{t,s}:=\bbE[X_t\,X_s]$. We have in mind
the case of $X$ being a Brownian motion or a Brownian bridge. In
this section we consider two processes $Y$ and $Z$, both defined by
linear transformations of $X$, and we write an absolute continuity
formula between the laws of $Y$ and $Z$.

For all $h$ in the space $M([0,1])$ of
all signed measures with finite total variation on $[0,1]$ we set:
\[
Q:M([0,1])\mapsto C([0,1]), \qquad Q\gl(t) \, := \, \int_0^1 q_{t,s} \, \gl(ds), \quad
t\in[0,1].
\]
We denote by $\langle \cdot,\cdot \rangle:C([0,1])\times M([0,1])\mapsto\bbR$ the canonical pairing,
\[
\langle h,\mu \rangle \, := \, \int_0^1 h_t \, \mu(dt).
\]
where a continuous function $k\in C([0,1])$ is identified with $k_t \,dt\in M([0,1])$.
We consider $\gl,\mu\in M([0,1])$ such that:
\begin{equation}\label{inde}
\langle Q \gl,\mu\rangle \, = \, 0, \qquad \langle Q\gl,\gl \rangle + \langle Q\mu,\mu\rangle\, = \, 1.
\end{equation}
We set for all $\go\in C([0,1])$:
\[
\gga(\go) \, := \, \int_0^1 \go_s \, \gl(ds), \qquad \Lambda_t \, := \, Q\gl(t), \ t\in[0,1],
\qquad I \, := \, \langle Q\gl,\gl \rangle,
\]
\[
a(\go) \, := \, \int_0^1 \go_s \, \mu(ds), \qquad M_t \, := \, Q\mu(t), \ t\in[0,1],
\qquad 1-I \, = \, \langle Q\mu,\mu \rangle,
\]
and we notice that $\gga(X)\sim N(0,I)$, $a(X)\sim N(0,1-I)$ and $\{\gga(X),a(X)\}$ are
independent by \eqref{inde}. We fix a constant $\kappa\in \bbR$ and
if $I<1$ we define the continuous processes
\[
Y_t\, := \, X_t+(\Lambda_t+M_t)\left(\kappa-a(X)-\gga(X) \right), \quad t\in[0,1],
\]
\[
Z_t\, := \, X_t+\frac 1{1-I}\, M_t\left(\kappa-a(X)-\gga(X) \right), \quad t\in[0,1].
\]
\begin{lemma}\label{abco} Suppose that $I<1$. Then
for all bounded Borel $\Phi:C([0,1])\mapsto\bbR$:
\begin{equation}\label{stga}
\bbE\left[\Phi(Y)\right] \, = \, \bbE\left[\Phi(Z) \, \rho(Z) \right],
\end{equation}
where for all $\go\in C([0,1])$:
\[
\rho(\go) \, := \, \frac 1{\sqrt{1-I}} \,
\exp\left(-\frac 12 \,\frac{1}{1-I} \,
\left(\gga(\go)-\kappa\right)^2+\, \frac12 \, \kappa^2 \right).
\]
\end{lemma}
\noindent
We postpone the proof of Lemma \ref{abco}
to section \ref{secproof}.

\section{The Brownian meander}

In this section we prove Theorem \ref{1.1}.
In the notation of section \ref{secabco}, we consider
$X=(B_t,t\in[0,1])$, standard Brownian motion.
It is easy to see that for all $t\in[0,1]$:
\[
\bbE\left[B_t \int_0^1 B_r \, dr\right] = \frac{t\,(2-t)}2,
\qquad \bbE\left[\left(\int_0^1 B_r \, dr\right)^2\right] = \frac 13.
\]
Therefore, it is standard that for all $c\in\bbR$, $B$ conditioned
to $\int_0^1 B=c$ is equal in law to the process:
\[
B^c_t \, := \, B_t \, + \, \frac 32\, t\, (2-t) \left(c-\int_0^1 B\right),
\qquad t\in[0,1].
\]
\begin{lemma}\label{l41}
Let $c\in\bbR$.
For all bounded Borel $\Phi:C([0,1])\mapsto\bbR$:
\[
\bbE\left[\Phi(B) \, \left| \, \int_0^1 B=c\right. \right] \, = \,
\bbE\left[\Phi(B^c)\right] \, = \, \bbE\left[\Phi\left(S\right)\, \rho(S)\right],
\]
where
\[
S_t \, := \, \left\{ \begin{array}{ll}
{\displaystyle
B_t, \qquad \qquad \qquad \qquad \qquad \qquad \qquad \quad t\in[0,1/2]
}
\\ \\
{\displaystyle
B_t+ \left(12\, t\, (2-t)-9\right)\left(c-\int_0^1 B\right), \quad t\in[1/2,1]}
\end{array} \right.
\]
\[
\rho(\go) \, := \, \sqrt 8 \, \exp\left(
-12\left(\int_0^\frac 12 \left(\go_r+\go_\frac 12\right) dr
- c\right)^2 +\frac 32 \, c^2 \right), \quad \go\in C([0,1]).
\]
\end{lemma}
\noindent{\bf Proof}. We are
going to show that we are in the setting of Lemma \ref{abco} with
$X=B$, $Y=B^c$ and $Z=S$.
We denote the Dirac mass at $\theta$ by $\delta_\theta$.
In the notation of section \ref{secabco}, we consider:
\[
\gl(dt) \, := \, \sqrt3 \left( 1_{[0,\frac 12]}(t) \, dt + \frac 12 \, \delta_{\frac 12}(dt)\right),
\qquad \mu(dt) \, := \, \sqrt3 \left( 1_{[\frac 12,1]}(t) \, dt - \frac 12 \, \delta_{\frac 12}(dt)\right),
\]
and $\kappa:=\sqrt 3 \, c$. Then:
\[
\gga(\go) \, = \, \sqrt 3\int_0^{\frac 12} \left(\go_r+\go_{\frac 12}\right) dr, \qquad
a(\go) \, = \, \sqrt 3\int_{\frac 12}^1 \left(\go_r-\go_{\frac 12}\right) dr,
\]
\[
\gga(\go)+a(\go) \, = \, \sqrt 3 \int_0^1 \go_r \, dr,
\qquad I \, = \, 3\int_0^{\frac 12} (1-r)^2 \, dr \, = \, \frac78.
\]
\[
\Lambda_t \, = \, \left\{ \begin{array}{ll}
{\displaystyle
\sqrt3 \ t \left(1-\frac {t}2\right), \quad t\in[0,1/2]
}
\\ \\
{\displaystyle \frac{3\sqrt 3}8,
\qquad \qquad t\in[1/2,1].}
\end{array} \right.
\
M_t \, = \, \left\{ \begin{array}{ll}
{\displaystyle
0, \ t\in[0,1/2]
}
\\ \\
{\displaystyle \sqrt3 \ t \left(1-\frac {t}2\right)-\frac{3\sqrt 3}8,
\quad t\in[1/2,1].}
\end{array} \right.
\]
Tedious but straightforward computations show that with these definitions we have
$X=B$, $Y=B^c$ and $Z=S$ in the notation of Lemma \ref{abco} and (\ref{inde}) holds true.
Then the thesis of Lemma \ref{l41} follows from
Lemma \ref{abco}. \qed

\medskip
\noindent{\bf Proof of Theorem \ref{1.1}}. Recall that $m$ is equal
in law to $B$ conditioned to be non-negative (see \cite{dim} and
\eqref{scaling} below). We want to condition $B$ first to be
non-negative and then to have a fixed time average. It turns out
that Lemma \ref{l41} allows to compute the resulting law by
inverting the two operations: first we condition $B$ to have a fixed
average, then we use the absolute continuity between the law of
$B^c$ and the law of $S$ and finally we condition $S$ to be
non-begative.

We set $K_\gep:=\{\go\in C([0,1]): \go\geq -\gep\}$,
$\gep\geq 0$. We recall that $B$ conditioned on $K_\gep$ tends
in law to $m$ as $\gep\to 0$, more generally for all $s>0$ and
bounded continuous $\Phi:C([0,s])\mapsto\bbR$, by the Brownian scaling:
\begin{equation}\label{scaling}
\lim_{\gep\to 0} \, \bbE\left[ \Phi(B_t, \, t\in[0,s]) \, \Big| \,
B_t\geq -\gep, \, \forall \, t\in[0,s] \right] \, = \,
\bbE\left[ \Phi\left( {\sqrt s}\, m_{t/s}, \, t\in[0,s] \right) \right],
\end{equation}
and this is a result of \cite{dim}.
By the reflection principle, for all $s>0$:
\begin{equation}\label{reflection}
\bbP(B_t\geq -\gep, \ \forall t\in[0,s]) \, = \, \bbP(|B_s|\leq \gep) \, \sim \,
\sqrt\frac2{\pi\, s} \, \gep, \qquad \gep\to 0.
\end{equation}
In particular for all bounded $f\in C(\bbR)$
\[
\bbE\left[\Phi(m) \, f(\langle m,1\rangle)\right]
\, = \, \lim_{\gep\to 0}\,
\sqrt\frac\pi2 \, \frac 1\gep \, \bbE\left[\Phi(B) \, 1_{K_\gep}(B) \, f(\langle B,1\rangle)\right].
\]
We want
to compute the limit of $\frac 1\gep \, \bbE\left[\Phi(B^c) \, 1_{K_\gep}(B^c)\right]$
as $\gep\to 0$.
Notice that $S$, defined in Lemma \ref{l41}, is equal to $B$ on
$[0,1/2]$. Therefore, by (\ref{scaling}) and (\ref{reflection}) with $s=1/2$:
\[
\sqrt\frac\pi2 \, \frac 1\gep \, \bbE\left[\Phi(B^c) \, 1_{K_\gep}(B^c)\right] \, \to \,
\sqrt 2 \, \bbE\left[\Phi(U^c)\, \rho(U^c) \, 1_{K_0}(U^c)\right].
\]
Comparing the last two formulae for all $f\in C(\bbR)$ with compact support:
\begin{align*}
& \sqrt\frac\pi2 \, \frac 1\gep \, \bbE\left[\Phi(B) \, 1_{K_\gep}(B) \, f(\langle B,1\rangle)\right]
\, = \, \int_\bbR \sqrt\frac\pi2 \, \frac 1\gep \, \bbE\left[\Phi(B^c) \, 1_{K_\gep}(B^c) \right]
f(c) \, N(0,1/3)(dc)
\\ & \to \, \int_0^\infty \sqrt\frac{24}\pi \,
\bbE\left[\Phi(U^c)\,
e^{-12\left(\int_0^{1/2}(U_r^c+U_{1/2}^c)\,
dr-c\right)^2}\, 1_{K_0}(U^c)\right]  \, f(c) \, dc
 = \bbE\left[\Phi(m) \, f(\langle m,1\rangle)\right]
\end{align*}
and (\ref{e1.1}) is proven. \qed

\section{The Brownian excursion}
\label{exx}

Let $(e_t,t\in[0,1])$ be the normalized Brownian excursion, see \cite{reyo}, and
$(\beta_t,t\in[0,1])$ a Brownian bridge between $0$ and $0$.
Let $\{m,\hat m,b\}$ be a triple of processes such that:
\begin{enumerate}
\item $m$ and $\hat m$ are independent copies of a Brownian meander on
$[0,1]$
\item conditionally on $\{m,\hat m\}$, $b$ is a Brownian bridge on
$[1/3,2/3]$ from $\frac 1{\sqrt 3}\, m_1$ to $\frac 1{\sqrt 3}\, \hat m_1$
\end{enumerate}
We introduce the continuous processes:
\begin{align*}
& v_t \, := \, \left\{ \begin{array}{ll}
{\displaystyle
\frac 1{\sqrt 3}\, m_{3t}, \qquad \quad t\in[0,1/3]
}
\\ \\
{\displaystyle
b_t, \qquad \qquad t\in[1/3,2/3],}
\\ \\
{\displaystyle
\frac 1{\sqrt 3}\, \hat m_{1-3t}, \quad t\in[2/3,1],}
\end{array} \right.
\\
& V^c_t \, := \,
\left\{ \begin{array}{ll}
{\displaystyle
v_t, \qquad \quad t\in[0,1/3]\cup[2/3,1]
}
\\ \\
{\displaystyle v_t
+ 18\left(9\, t\, (1-t)-2\right)\left(c-\int_0^1 v\right), \quad t\in[1/3,2/3].}
\end{array} \right.
\end{align*}
Notice that $\int_0^1 V^c_t \, dt = c$.
We set for all $\go\in C([0,1])$:
\[
\rho^c(\go) \,:= \, \exp\left\{
-162\left(\int_0^\frac 13 \left(\go_r+\go_{1-r}\right) dr + \frac{\go_\frac 13
+\go_\frac 23}6- c\right)^2 - \frac 32 (\go_{\frac 23}-\go_{\frac 13})^2 \right\}.
\]
\begin{theorem}\label{4.1}
For all bounded Borel $\Phi:C([0,1])\mapsto\bbR$ and $f:\bbR\mapsto\bbR$:
\begin{equation}\label{e4.1}
\bbE\left[\Phi(e) \, f(\langle e,1\rangle)\right] \, = \,
\int_0^\infty 27\sqrt{\frac6{\pi^3}} \, \bbE\left[\Phi\left(V^c\right)\, \rho^c(V^c) \,
\, 1_{K_0}(V^c)\right] \, f(c) \, dc
\end{equation}
\end{theorem}
\begin{cor}\label{cor3}
If we set
\[
p_{\langle e,1\rangle}(c) \, =: \, 27\sqrt{\frac6{\pi^3}} \,
\bbE\left[\rho^c(V^c) \, 1_{\{V^c_t\geq 0, \ \forall t\in[0,1]\}}
\right].
\]
then $p_{\langle e,1\rangle}$ is the density of $\langle e,1\rangle$
on $[0,\infty)$, i.e.
\[
\bbP(\langle e,1\rangle\in dc)\, = \, p_{\langle e,1\rangle}(c)
\,1_{\{c\geq 0\}} \, dc.
\]
Moreover $p_{\langle e,1\rangle}$ is continuous on $[0,\infty)$.
\end{cor}
\noindent Notice that a.s. $V^c_t\geq 0$ for all
$t\in[0,1/3]\cup[2/3,1]$, since a.s. $m\geq 0$: therefore a.s.
\[
\{V^c_t\geq 0, \ \forall t\in[0,1]\} \, = \, \{V^c_t\geq 0, \
\forall t\in[1/3,2/3]\}.
\]
The probability of this event is positive for all $c>0$ while it is
0 for $c=0$, since $\int_0^1 V^0_t \, dt = 0$. In particular
$p_{\langle e,1\rangle}(0)=0$. Moreover for all $c> 0$ we can define
the probability measure on $C([0,1])$ given by
\[
\bbE\left[\Phi(e) \, | \, \langle e,1\rangle =c \right] \, := \,
\frac 1{\cZ_c} \ \bbE\left[\Phi\left(V^c\right)\, \rho^c(V^c) \, \,
1_{\{V^c_t\geq 0, \ \forall t\in[0,1]\}} \right]
\]
where $\Phi:C([0,1])\mapsto\bbR$ is bounded Borel and $\cZ_c>0$ is a
normalization factor.
\begin{cor}\label{cor4}
$(\bbP\left[e\in\cdot \, | \, \langle e,1\rangle =c \right], c>0)$
is a regular conditional distribution of $e$ given $\langle
e,1\rangle$, i.e.
\[
\bbP(e\in\cdot\, , \langle e,1\rangle\in dc ) \, = \,
\bbP\left[e\in\cdot \, | \, \langle e,1\rangle =c \right] \
p_{\langle e,1\rangle}(c) \,1_{\{c> 0\}} \, dc.
\]
Moreover $(0,\infty)\ni c\mapsto \bbP\left[e\in\cdot \, | \, \langle
e,1\rangle =c \right]$ is continuous in the weak topology.
\end{cor}
%

\medskip\noindent
Notice that for all $t\in[0,1]$:
\[
\bbE\left[\beta_t \int_0^1 \beta_r \, dr\right] = \frac{t(1-t)}2,
\qquad \bbE\left[\left(\int_0^1 \beta_r \, dr\right)^2\right] = \frac 1{12}.
\]
Therefore, for all $c\in\bbR$, $\beta$ conditioned
to $\int_0^1 \beta=c$ is equal in law to the process:
\[
\beta^c_t \, := \, \beta_t \, + \, 6\, t\, (1-t) \left(c-\int_0^1 \beta\right),
\qquad t\in[0,1].
\]
\begin{lemma}\label{l44}
Let $c\in\bbR$.
For all bounded Borel $\Phi:C([0,1])\mapsto\bbR$:
\[
\bbE\left[\Phi(\beta) \, \left| \, \int_0^1\beta=c \right. \right]  \, = \,
\bbE\left[\Phi(\beta^c)\right]
\, = \, \bbE\left[\Phi\left(\Gamma^\beta\right)\, \rho_1\left(\Gamma^\beta\right)\right]
\]
where for all $\go\in C([0,1])$
\[
\Gamma^\go_t \, = \, \left\{ \begin{array}{ll}
{\displaystyle
\go_t, \qquad \qquad \qquad \qquad \qquad \qquad \qquad  \quad t\in[0,1/3]\cup[2/3,1]
}
\\ \\
{\displaystyle
\go_t+ 18\left(9\, t\, (1-t)-2\right)\left(c-\int_0^1 \go\right), \quad t\in[1/3,2/3]}
\end{array} \right.
\]
\[
\rho_1(\go) \,:= \, \sqrt{27} \, \exp\left(
-162\left(\int_0^\frac 13 \left(\go_r+\go_{1-r}\right) dr + \frac{\go_\frac 13
+\go_\frac 23}6- c\right)^2 +6\, c^2 \right).
\]
\end{lemma}
\noindent
{\bf Proof}. Similarly to the proof of Lemma \ref{l41}, we are going to show
that we are in the situation of Lemma \ref{abco} with $X=\beta$, $Y=\beta^c$ and
$Z=\Gamma^\beta$.
In the notation of Lemma \ref{abco}, we consider
\[
\gl(dt) \, := \, \sqrt{12} \left( 1_{[0,\frac 13]\cup[\frac 23,1]}(t) \, dt +
\frac{\delta_{\frac 13}(dt)+\delta_{\frac 23}(dt)}6 \right),
\]
\[
\mu(dt) \, := \, \sqrt{12} \left( 1_{[\frac 13,\frac 23]}(t) \, dt -
\frac{\delta_{\frac 13}(dt)+\delta_{\frac 23}(dt)}6 \right),
\]
and $\kappa:=\sqrt{12} \, c$. Then:
\[
\gga(\beta) \, = \, \sqrt{12}\int_0^{\frac 13} \left(\beta_r+\frac 12 \, \beta_{\frac 13}\right) dr
+ \sqrt{12}\int_{\frac 23}^1 \left(\beta_r+\frac 12 \, \beta_{\frac 23}\right) dr,
\quad I \, = \, \frac{26}{27}
\]
\[
a(\beta) \, = \, \sqrt{12}\int_{\frac 13}^{\frac 23} \left(\beta_r-\frac 12 \, \beta_{\frac 13}
- \frac 12 \, \beta_{\frac 23}\right) dr,
\]
\[
\Lambda_t \, = \ 1_{[0,\frac 13]\cup[\frac 23,1]}(t) \ \sqrt3 \, t(1-t) \
+ 1_{(\frac 13,\frac 23)}(t) \ \frac{2\sqrt{3}}9,
\quad
M_t \, = \ 1_{[\frac 13,\frac 23]}(t) \ \sqrt3 \, t(1-t).
\]
Again the thesis follows by direct computations and from
Lemma \ref{abco}. \qed

\medskip\noindent
{\bf Proof of Theorem \ref{4.1}}. We follow the proof of
Theorem \ref{1.1}.
Define $\{B,b,\hat B\}$, processes such that:
\begin{enumerate}
\item $B$ and $\hat B$ are independent copies of a standard Brownian motion over
$[0,1/3]$
\item conditionally on $\{B,\hat B\}$, $b$ is a Brownian bridge over $[1/3,2/3]$ from
$B_{1/3}$ to $\hat B_{1/3}$.
\end{enumerate}
We set:
\[
r_t \, := \, \left\{ \begin{array}{ll}
B_t \qquad \qquad t\in[0,1/3]
\\ \\
b_t \qquad \qquad t\in[1/3,2/3]
\\ \\
\hat B_{1-t} \qquad \qquad t\in[2/3,1].
\end{array} \right.
\]
Moreover we set, denoting the density of $N(0,t)(dy)$ by $p_t(y)$:
\[
\rho_2(\go) \, := \,
\frac{p_\frac 13(\go_\frac  23-\go_\frac  13)}{p_1(0)}
\, = \,  \sqrt{3} \,
\exp\left( - \frac 32 \, (\go_{\frac 23}-\go_{\frac 13})^2\right),
\qquad \go\in C([0,1]).
\]
By the Markov property of $\beta$:
\[
\bbE \left[\Phi(r) \, \rho_2(r)\right] \, = \, \bbE[\Phi(\beta)].
\]
Then, recalling the definition of $\rho^c$ above, by Lemma
\ref{abco} and Lemma \ref{l44}:
\[
\bbE \left[\Phi(\beta^c)\right] \, = \,
\bbE\left[\Phi\left(\Gamma^\beta\right)\, \rho_1\left(\Gamma^\beta\right)\right]
\, = \, \bbE[\Phi\left(\Gamma^r\right)\, \rho_1\left(\Gamma^r\right)
\, \rho_2(\Gamma^r)] \, = \, 9 \
\bbE[\Phi(\Gamma^r) \, \rho^c(\Gamma^r)] \, e^{6c^2}.
\]
We recall now that $\bbP(\beta\in K_\gep)=1-\exp(-2\, \gep^2) \sim
2\, \gep^2$ as $\gep\to0$, where $K_\gep=\{\go\in C([0,1]): \go\geq
-\gep\}$. We want to compute the limit of $\frac 1{2\, \gep^2} \,
\bbE\left[\Phi(\beta^c) \, 1_{K_\gep}(\beta^c)\right]$ as $\gep\to
0$. On the other hand $\bbP(B_t\geq -\gep, \forall t\in[0,1/3])\sim
\sqrt{\frac6\pi}\, \gep$ by (\ref{reflection}).  Then by
(\ref{scaling}) and (\ref{reflection})
\[
\frac 1{2\, \gep^2} \, \bbE\left[\Phi(\beta^c) \, 1_{K_\gep}(\beta^c)\right] \, \to \,
\frac{27}\pi \, \bbE\left[\Phi\left(V^c\right)\, \rho^c(V^c) \,
\, 1_{K_0}(V^c)\right] e^{6c^2}, \quad \frac{27}\pi \, = \, \frac 12 \, \sqrt{\frac6\pi}\, \sqrt{\frac6\pi}\, 9.
\]
On the other hand, $\beta$ conditioned on $K_\gep$ tends
in law to the normalized Brownian excursion
$(e_t,t\in[0,1])$, as proven in \cite{dim}. Then we have for all bounded $f\in C(\bbR)$:
\[
\frac 1{2\, \gep^2} \, \bbE\left[\Phi(\beta) \, 1_{K_\gep}(\beta) \, f(\langle \beta,1\rangle)\right]
\, \to \, \bbE\left[\Phi(e) \, f(\langle e,1\rangle)\right]
\]
Comparing the two formulae for all $f\in C(\bbR)$ with compact support:
\begin{align*}
& \frac 1{2\, \gep^2} \, \bbE\left[\Phi(\beta) \, 1_{K_\gep}(\beta) \, f(\langle \beta,1\rangle)\right]
\, = \, \int_\bbR \frac 1{2\, \gep^2} \, \bbE\left[\Phi(\beta^c) \, 1_{K_\gep}(\beta^c) \right]
f(c) \, N(0,1/12)(dc)
\\ & \to \, \int_0^\infty 27\sqrt{\frac6{\pi^3}} \, \bbE\left[\Phi\left(V^c\right)\, \rho^c(V^c) \,
\, 1_{K_0}(V^c)\right] \, f(c) \, dc \, = \, \bbE\left[\Phi(e) \, f(\langle e,1\rangle)\right]
\end{align*}
and (\ref{e4.1}) is proven. \qed

\section{Proof of Proposition \ref{abco}.}
\label{secproof}

The thesis follows if we show that the Laplace
transforms of the two probability measures in (\ref{stga}) are equal.
Notice that $Y$ is a Gaussian process with mean $\kappa\, (\Lambda+M)$ and covariance function:
\[
q^Y_{t,s} \, = \, \bbE\left[\left(Y_t-\kappa\,(\Lambda_t+M_t)\right)\left(Y_s-\kappa\,(\Lambda_s+M_s)\right)\right] \, = \,
q_{t,s} \, - \, (\Lambda_t+M_t)\,(\Lambda_s+M_s),
\]
for $t,s\in[0,1]$. Therefore, setting for all $h\in C([0,1])$: $Q_Yh(t):=\int_0^1 q^Y_{t,s} \, h_s \, ds$,
$t\in[0,1]$, the Laplace transform of the law of $Y$ is:
\[
\bbE\, \left[e^{\langle Y,h\rangle}\right] \, = \, e^{\kappa\langle h,\Lambda+M\rangle+
\frac12 \, \langle Q_Y h,h\rangle}.
\]
Recall now the following version of the Cameron-Martin Theorem: for all $h\in M([0,1])$
\[
{\mathbb E}\left[ \Phi(X) \, e^{\langle X,h\rangle} \right] \, = \,
e^{\frac 12\langle Qh,h\rangle} \, {\mathbb E} [\Phi(X+Qh)].
\]
Notice that $\gga(Z)=\gga(X)$, by (\ref{inde}). Therefore $\rho(Z)=\rho(X)$.
We obtain, setting $\overline h:=h-\frac1{1-I}\langle M,h\rangle(\gl+\mu)$:
\begin{align*}
& \bbE\left[e^{\langle Z,h\rangle} \, \rho(Z) \right] \, = \, e^{\frac\kappa{1-I}\langle M,h\rangle} \,
\bbE\left[e^{\langle X,\overline h\rangle}\, \rho(X) \right]
\, = \, e^{\frac\kappa{1-I}\langle M,h\rangle+\frac 12\langle Q\overline h,\overline h\rangle} \, {\mathbb E}
\left[\rho\left(X+Q\overline h\right)\right] =
\\ & = \, e^{\frac\kappa{1-I}\langle M,h\rangle+\frac 12\langle Q\overline h,\overline h\rangle} \,
\frac 1{\sqrt{1-I}} \, \bbE\left[
e^{-\frac 12 \, \frac 1{1-I} \, \left(\gga(X)+\langle \overline h,\Lambda\rangle-\kappa\right)^2+
\frac12 \, \kappa^2}\right].
\end{align*}
By the following standard
Gaussian formula for $\ga\sim N(0,\sigma^2)$, $\sigma\geq 0$ and $c\in\bbR$:
\[
\bbE\left[e^{-\frac 12 \, (\ga+c)^2}\right]
\, = \, \frac 1{\sqrt{1+\sigma^2}} \, e^{-\frac 12
\, \frac{c^2}{1+\sigma^2}},
\]
we have now for $\gga(X)\sim N(0,I)$:
\begin{align*}
& \bbE\left[
e^{-\frac 12 \, \frac 1{1-I} \, \left(\gga(X)+\langle \overline h,\Lambda\rangle-\kappa\right)^2}
\right] \, = \, \frac1{\sqrt{1+\frac{I}{1-I}}} \, e^{-\frac 12 \, \frac 1{1-I} \,
\frac1{1+\frac{I}{1-I}}\left(\langle \overline h,\Lambda\rangle-\kappa\right)^2}
\, = \, {\sqrt{1-I}} \, e^{-\frac 12 \, \left(\langle \overline h,\Lambda\rangle-\kappa\right)^2}
\end{align*}
Therefore, recalling the definition of $\overline
h:=h-\frac1{1-I}\langle M,h\rangle(\gl+\mu)$, we obtain after some
trivial computation:
\begin{align*} &
\log \bbE\left[e^{\langle Z,h\rangle} \, \rho(Z) \right] \, = \,
\frac\kappa{1-I}\langle M,h\rangle+\frac 12\langle Q\overline h,\overline h\rangle
-\frac 12 \, \left(\langle \overline h,\Lambda\rangle-\kappa\right)^2 + \frac 12 \, \kappa^2
\\ & = \kappa\langle \Lambda+M,h\rangle + \frac 12 \langle Qh,h\rangle
- \langle \Lambda+M,h\rangle^2
\, = \, \kappa\langle h,\Lambda+M\rangle+\frac 12\langle Q_Yh,h\rangle. \qed
\end{align*}

\end{document}